\documentclass[12pt]{article}

\usepackage{amsthm,amsmath,amssymb}
\usepackage{graphicx}
\usepackage[left=25mm, right=25mm,
            top=30mm, bottom=30mm]{geometry}

\theoremstyle{plain}
\newtheorem{theorem}{Theorem}
\newtheorem{lemma}[theorem]{Lemma}

\theoremstyle{definition}

\theoremstyle{remark}
\newtheorem{remark}[theorem]{Remark}

\title{\bf The Page-R\'enyi parking process}

\author{Lucas Gerin\thanks{Supported by grant \textsc{Anr Graal}.}\\
\small CMAP, \'Ecole Polytechnique\\[-0.8ex]
\small Route de Saclay,\\[-0.8ex] 
\small 91128 Palaiseau, France\\
\small\tt gerin@cmap.polytechnique.fr
}

\date{\small Mathematics Subject Classifications: 60C05,68R05}

\begin{document}

\maketitle

\begin{abstract}
In the Page parking (or packing) model on a discrete interval (also known as the \emph{discrete R\'enyi packing problem} or the \emph{unfriendly seating problem}), cars of length two successively park uniformly at random on pairs of adjacent places, until only isolated places remain.

We use a probabilistic construction of the Page parking to give a short proof of the (known) fact that the proportion of the interval occupied by cars goes to $1-e^{-2}$, when the length of the interval goes to infinity.
We also obtain some new consequences on both finite and infinite parkings.

\end{abstract}

\section{The Page parking}
\newcommand{\set}[1]{\left\{#1\right\}}
\newcommand{\eps}{\varepsilon}
\newcommand{\lf}{\lfloor}
\newcommand{\rf}{\rfloor}
\subsection{The model}
For $n\geq 2$, we consider a sequence of \emph{parking configurations} $x^t=\left(x^t_i\right)_{1\leq i\leq n}$ in $\set{0,1}^n$, given by the following construction. Initially the parking is empty: $x^0=0^n$. Given $x^t$ one draws uniformly at random (and independently from the past) a number $i$ in $\set{1,2,\dots, n-1}$ and, if possible, a car of size $2$ \emph{parks} at $i$:
\begin{itemize}
\item if $(x^t_i,x^t_{i+1})= (0,0)$ then $(x^{t+1}_i,x^{t+1}_{i+1})=(1,1)$ and $n-2$ other coordinates remain unchanged;
\item if $(x^t_i,x^t_{i+1})\neq (0,0)$ nothing happens.
\end{itemize}
After some random time $T_n$ (which is dominated by a coupon collector process with $n-1$ coupons, see Section \eqref{Sec:Coupons} below) 
parking is no longer possible, in the sense that in $x^{T_n}$ there are no adjacent coordinates  $(i,i+1)$ such that $x^{T_n}_i=x^{T_n}_{i+1}=0$. We set $X_{n}=x^{T_n}$ and $X_{n}(i)=x^{T_n}_i$. Below is an example where $X_n=(1,1,0,1,1,1,1,0,1,1,0,1,1)$.

\begin{center}
\includegraphics[width=12cm]{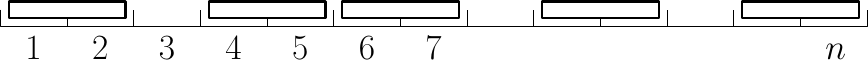}\\
\end{center}

We are mainly interested in the numbers $M_n$ of places occupied by a car:
$$
M_{n}=\mathrm{card}\set{1\leq i\leq n,\ X_{n}(i)=1}.
$$
We obviously have $2\lf n/3\rf \leq M_{n}\leq n$ and we expect $M_{n}/n$ to converge, at least in some sense. Page obtained the following law of large numbers.
\begin{theorem}[Page (1959)]\label{Th:LimitePage}
When $n\to +\infty$,
$$
\frac{M_{n}}{n} \stackrel{\text{prob.}}{\to} 1-e^{-2}=0.8646647\dots
$$
\end{theorem}
The fact that $\mathbb{E}[M_n]/n$ converges to $1-e^{-2}$ was in fact already observed by Flory \cite{Flory} in 1939.
The proof of Theorem \ref{Th:LimitePage} is essentially obtained on conditioning on the position $i$ of the first car. This gives the recursion identity
$$
M_{n}\stackrel{\text{(d)}}{=} M_{I-1}+M'_{n-I-1}+2,
$$
where $I$ is uniform in $\set{1,2,\dots,n-1}$ and $M_{I-1},M'_{n-I-1}$ are independent conditionally on $I$. This gives recursions for both sequences $\mathbb{E}[M_n]$ and $\mathbb{E}[M_n^2]$, which can be handled using generating functions. Theorem \ref{Th:LimitePage} can also be found in \cite{Friedman,Fan,Flajolet,Pinsky} with similar proofs.

The Page parking problem has a long story, it has been studied by many people and under different names. It is equivalent to the \emph{unfriendly seating} problem \cite{Friedman}, and sometimes also called the \emph{discrete R\'enyi Packing} model \cite{Hemmer}, more generally, it is a toy model for \emph{random deposition}.
We refer to \cite{Flory,Page,Chimie} for some interpretations of the model in polymer chemistry.
Besides, we mention that much is known when cars have size $\ell >2$, see \cite{KlassenRunnenburg,Flajolet,Pinsky}. 

The first aim of the present paper is to present a probabilistic (and new) proof of Theorem \ref{Th:LimitePage} and to study the asymptotic behavior of $T_n$ (the numbers of cars that have tried to park). We also prove a law of large numbers for the infinite Page parking and extend a formula due to Hemmer \cite{Hemmer} for a continuous-time version of the Page parking.

\subsection{The probabilistic construction}

An alternative way of defining $M_{n}$ is the following. Let $\xi=(\xi_i)_{1\leq i\leq n-1}$ be i.i.d. non-negative random variables with continuous distribution function $F$ (hereafter we will take $F(t)=1-e^{-t}$). By convention we set $\xi_0=\xi_n=+\infty$.

Then the order statistics $\left(\xi_{\sigma(1)}< \dots <\xi_{\sigma(n-1)}\right)$ give the order in which the cars park:
\begin{itemize}
\item at time $t=\xi_{\sigma(1)}$, the first car parks at $(\sigma(1),\sigma(1)+1)$,
\item at time $t=\xi_{\sigma(2)}$, the second one parks (if possible) at $(\sigma(2),\sigma(2)+1)$,
\item ...
\item at time $t=\xi_{\sigma(n-1)}$, a car parks (if possible) at $(\sigma(n-1),\sigma(n-1)+1)$ and the process is over.
\end{itemize}
It is easy to see that we obtain the same distribution, this was already observed by previous authors (see \cite{Runnenburg}).
Let us collect for further use some obvious features of this construction.
\begin{remark}
\begin{itemize}
\item The configuration $X_n$ only depends on the ordering of the $\xi_i$'s.
\item If $i$ is a local minimum of $\xi$ (\emph{i.e.}if  $\xi_{i-1}>\xi_i<\xi_{i+1}$) then places $i$ and $i+1$ are empty at time $t=\xi_i^-$ and then a car parks at time $\xi_i$ at $(i,i+1)$.
\item To define $X_n$ from the $\xi_i$'s, one can treat separately the intervals defined by two successive local minima.
\end{itemize}
\end{remark}
Here is a sample of $\xi$ (here $\xi_{\sigma(1)}=\xi_6$) and the corresponding configuration $X_{n}$:
\begin{center}
\includegraphics[width=12cm]{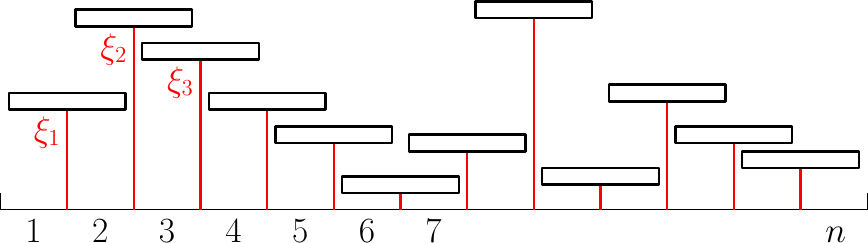}\\
\vspace{3mm}
$\downarrow$ \ \\
\vspace{3mm}
\includegraphics[width=12cm]{Config_finale.pdf}\\
\end{center}

It does not seem that this construction was used to its full extent, yet it gives a very simple way to characterize positions $i$ occupied by a car. We need a few definitions.

We say that there is a \emph{rise} of length $\ell$ at $i$ if $i\geq \ell+2$ and (recall $\xi_0=\xi_n=+\infty$)
$$
\xi_{i-\ell-1}> \xi_{i-\ell}< \xi_{i-\ell+1} < \xi_{i-\ell+2} <\dots < \xi_{i-1} 
$$
or if $i=\ell+1$ and
$$
\xi_{1}< \xi_{2} < \dots < \xi_{i-1}.
$$
There is a \emph{descent} of length $\ell$ at $i$ if $i<n-\ell$ and
$$
\xi_{i}> \xi_{i+1} >  \dots > \xi_{i+\ell-1} <\xi_{i+\ell} 
$$
or if $i=n-\ell$ and
$$
\xi_{i}> \xi_{i+1} >  \dots > \xi_{n-1}.
$$
Consistently we say that there is a rise (resp. a descent) of length $1$ at $i$ if $\xi_{i-2}>\xi_{i-1}$ (resp. if $\xi_{i}<\xi_{i+1}$), so that for every $i$ there is a rise and a descent at $i$.

Note that by construction for each $i,\ell,\ell'$ the two events $\set{\text{ rise of length $\ell$ at $i$ }}$ and $\set{\text{ descent of length $\ell'$ at $i$ }}$ are independent.

\begin{lemma}\label{Lem:RiseDescent}
There is no car at $i$ (\emph{i.e.} $X_{n}(i)=0$) if and only if there is a rise of even length at $i$ and a descent of even length at $i$.
\end{lemma}

Here is an example of a rise of length $2\ell=6$ at $i$ and a descent of length $2k=4$:
\begin{center}
\includegraphics[width=10cm]{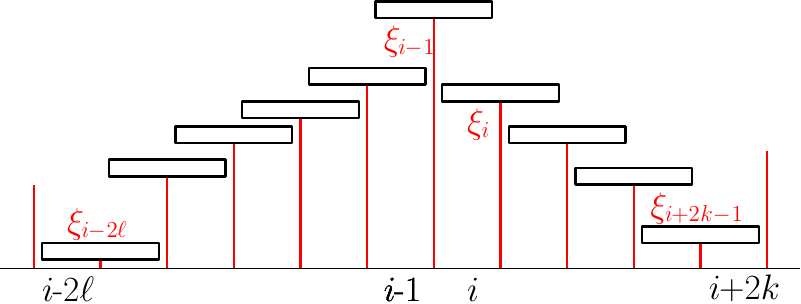}\\
\vspace{3mm}
$\downarrow$ \ \\
\vspace{3mm}
\includegraphics[width=10cm]{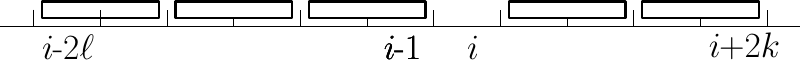}\\
\end{center}

\begin{proof}[Proof of Lemma \ref{Lem:RiseDescent}]
As already said, if $i$ is a local minimum then a car parks at $(i,i+1)$. Consider now the case where $i$ is inside a rise: $\xi_{r-1}>\xi_r<\xi_{r+1}<\xi_{r+2}<\dots <\xi_i <\xi_{i+1}$ for some $r$. Then cars successively park at $(r,r+1)$, $(r+2,r+3)$, $\dots$, until $(i-1,i)$ or $(i,i+1)$ (depending on the parity of $i-r$). Then $X_n(i)=1$, and the same applies in the case where $i$ is inside a descent.

Therefore the only case where $X_n(i)$ could be zero is if $i$ (or $i-1$) is a local maximum for $\xi$.
Define $m_i$ (resp. $m'_i$) as the closest local minimum of $\xi$ on the left of $i$ (resp. on the right). A rise begins at $m_i$ and a descent ends at $m'_i$ and $m_i=i-s$, $m'-i=s'-1$ where $s,s'$ are the lengths of these rise and descent.

Cars successively try to fill places of the rise $(m_i,m_{i}+1)$, $(m_{i}+2,m_{i}+3)$, $\dots$ from left to right and places of the descent $(m'_i,m'_{i}+1)$, $(m'_i-2,m'_{i}-1)$, $\dots$ from right to left. If only $s$ is odd (resp. only $s'$) then $i$ is occupied by a car of the rise (resp. descent). If both are odd a car parks at $(i-1,i)$ or $(i,i+1)$ depending on wether $\xi_{i-1}<\xi_i$ or not. If both are even the rightmost car of the rise parks at $(i-2,i-1)$ and the leftmost car of the descent parks at $(i+1,i+2)$, leaving $i$ unoccupied.

\end{proof}

\section{The infinite parking}
An interesting feature of the probabilistic construction is that it allows us to define the model $\left(X_\infty(i)\right)_{i\in\mathbb{Z}}$ on $\mathbb{Z}$, by considering a doubly-infinite sequence $\left(\xi_{i}\right)_{i\in \mathbb{Z}}$. 

We first set $X_\infty(i)=X_\infty(i+1)=1$ for every $i$ such that $\xi_i$ is a local minimum. Then we define $X_\infty(i)$ as before, using only the $\xi$'s between $m_i,m'_i$. Recall that $m_i$ (resp. $m'_i$) is the closest local minimum of $\xi$ on the left of $i$ (resp. on the right) and note that with probability one, for all $i$ one has $m'_i -m_i< +\infty$. Lemma \ref{Lem:RiseDescent} also holds for the infinite parking.
\begin{theorem}\label{Theo:Dimer}
In the infinite Page parking:
\begin{enumerate}
\item For every $i\in\mathbb{Z}$,
$$
\mathbb{P}(X_{\infty}(i)=1)=1-e^{-2}.
$$
\item When $n\to +\infty$,
$$
\frac{X_{\infty}(1)+\dots +X_{\infty}(n)}{n} \stackrel{\text{prob.}}{\to} 1-e^{-2}.
$$
\end{enumerate}
\end{theorem}
\begin{proof} We first prove the first statement.\\
By construction the rise at $i$ and the descent at $i$ are independent. By Lemma \ref{Lem:RiseDescent} we have
\begin{align*}
\mathbb{P}(X_{n}(i)=0)&=\mathbb{P}(\text{ even rise at }i\ ) \mathbb{P}(\text{ even descent at }i\ )\\
&=\mathbb{P}(\text{ even descent at }i\ )^2\\
&= \left(\sum_{\ell\geq 1} \mathbb{P}\left(\xi_{i}> \xi_{i+1} >  \dots > \xi_{i+2\ell-1} <\xi_{i+2\ell}  \right)\right)^2\\
&= \left(\sum_{\ell\geq 1} \mathbb{P}\left(\sigma_1> \sigma_2 >  \dots > \sigma_{2\ell} <\sigma_{2\ell+1}  \right)\right)^2,
\end{align*}
where $\sigma$ is a uniform permutation of $2\ell+1$ elements. We get
\begin{align*}
\mathbb{P}(X_{n}(i)=0)&= \left(\sum_{\ell\geq 1} \mathbb{P}\left(\sigma_1> \sigma_2 >  \dots > \sigma_{2\ell}\right)-\mathbb{P}\left(\sigma_1> \sigma_2 >  \dots > \sigma_{2\ell}>\sigma_{2\ell+1}  \right)\right)^2\\
&=\left(\sum_{\ell\geq 1} \frac{1}{(2\ell)!}-\frac{1}{(2\ell+1)!}\right)^2 = \left(\sum_{\ell\geq 0} \frac{1}{(2\ell)!}-\frac{1}{(2\ell+1)!}\right)^2= (1/e)^2.
\end{align*}
We now turn to the law of large numbers. We will prove that there exists a constant $c>0$ such that for every $i,j$
\begin{equation}\label{eq:cov}
\left|\mathrm{cov}(X_{\infty}(i),X_{\infty}(j))\right|\leq \frac{c}{\lf|j-i|/3\rf!}.
\end{equation}
The convergence in probability will follow since \eqref{eq:cov} implies that 
$$
\mathrm{Var}(X_{\infty}(1)+\dots +X_{\infty}(n))=\mathcal{O}(n).
$$
which gives the desired convergence.

Let us prove \eqref{eq:cov}. First note that $\mathrm{cov}(X_{\infty}(i),X_{\infty}(j))=\mathbb{P}(X_{\infty}(i)=0,X_{\infty}(j)=0)-e^{-4}$. By translation-invariance we can fix $i=1$ $j=p$ and take $p$ larger than, say, $10$,
\begin{align*}
0&\leq \mathbb{P}(X_{\infty}(1)=0,X_{\infty}(p)=0)\\
&-\mathbb{P}
\begin{pmatrix}
\text{ even rise at }1,\text{ even descent smaller than $\lf p/3\rf$ at }1, \\
\text{ even rise smaller than $p/3$ at }p,\text{ even descent at }p.
\end{pmatrix}\\
&\leq \mathbb{P}(\text{ even descent larger than $\lf p/3\rf$ at }1\ \cup\text{ even rise larger than $\lf p/3\rf$ at }p)
\end{align*}
and the latter probability is less than $2/\lf p/3\rf!$. Now, 
\begin{multline*}
\mathbb{P}
\begin{pmatrix}
\text{ even rise at }1,\text{ even descent smaller than $\lf p/3\rf$ at }1, \\
\text{ even rise smaller than $\lf p/3\rf$ at }p,\text{ even descent at }p.
\end{pmatrix}\\
= e^{-1}\times\left(\sum_{\ell:\ 2\ell\leq \lf p/3\rf}\frac{2\ell}{(2\ell +1)!}\right)^2 \times e^{-1}
= e^{-4}+\mathcal{O}\left(\frac{1}{(\lf p/3\rf)!} \right).
\end{multline*}
This proves \eqref{eq:cov} for $|j-i|\geq 10$, by taking a possibly larger $c$ we obtain the proof for every $i,j$.
\end{proof}

Note that the probabilistic construction also allows us to compute explicitly some short-range correlations. Clearly, by construction, if $X_{\infty}(1)=0$ then $X_{\infty}(2)=X_{\infty}(3)=1$, so the simplest non-trivial computation is
\begin{align*}
\mathbb{P}(X_{\infty}(1)=0,X_{\infty}(4)=0)&= \mathbb{P}(X_{\infty}(1)=0,X_{\infty}(2)=1,X_{\infty}(3)=1,X_{\infty}(4)=0)
\\&= \mathbb{P}(\text{ even rise at }1\ )\mathbb{P}\left(\xi_{1}> \xi_{2}<\xi_3\right) \mathbb{P}(\text{ even descent at }4\ )\\
&=e^{-2}/3.
\end{align*}
As one could expect, this is greater than $\mathbb{P}(X_{\infty}(1)=0)\times \mathbb{P}(X_{\infty}(4)=0)=e^{-4}$.

\subsection{Evolution of the continuous-time process}\label{Sec:density}
We now consider the process given by the time arrivals of cars. As above $X_{\infty}(i)$ is the indicator that there is eventually a car at $i$. We define the continuous-time process $\left(X^t_{\infty}\right)_{t\geq 0}$ with values in $\set{0,1}^{\mathbb{Z}}$ by
$$
X^t_{\infty}(i)=
\begin{cases}
1& \text{ if }X_{\infty}(i)=1\text{ and }\tau_i \leq t,\\
0& \text{ otherwise.}
\end{cases}
$$
Here, $\tau_i = \xi_{i-1}$ if the car parked at $i$ is parked at $(i-1,i)$ and $\tau_i = \xi_{i}$ if this car is at $(i,i+1)$. Then  $\tau_i$ is indeed the time arrival of the corresponding car, we set $\tau_i=+\infty$ if there is no car at $i$.

Recall that $F$ denotes the distribution function of the $\xi_i$'s. If $F(t)=1-e^{-t}$ then $(X^t_\infty)_{t\geq 0}$ defines a homogeneous Markov process and in this case the distribution function of $\tau_i$ was identified by Hemmer (\cite{Hemmer} eq. (19)). Here we generalize his result to any $F$.

\begin{theorem}[Evolution of the density of cars]\label{Th:EvolutionDensite}\ \\
Let $\tau_i$ be the arrival time of the car $i$,
$$
\mathbb{E}[X^t_{\infty}(i)]=\mathbb{P}(\tau_i\leq t)=1-e^{-2F(t)}.
$$
\end{theorem}
\noindent Of course we recover $\mathbb{P}(X_{\infty} =1)=\mathbb{P}(\tau_i< +\infty)= \lim_{t\to +\infty} 1-e^{-2F(t)}= 1-e^{-2}$. 
\begin{proof}
By translation-invariance we assume $i=0$.
Lemma \ref{Lem:RiseDescent} gives that $\tau_i\leq t$ if and only if
\begin{itemize}
\item there is an odd rise at $0$ and $\xi_{-1}\leq t$, 
\item or there is an odd descent at $0$ and $\xi_{0}\leq t$.
\end{itemize}
These two events being independent we have
$$
\mathbb{P}(\tau_i \leq t)=2f(t) -f(t)^2=f(t)\left(2-f(t)\right)
$$
where
$$
f(t)=\mathbb{P}(\xi_{0}\leq t; \text{ odd descent at $0$ }).
$$
Now,
\begin{multline*}
f(t)=  \mathbb{P}\left(t\geq \xi_{0}<\xi_1\right) + \sum_{k\geq 1} \mathbb{P}\left(t\geq \xi_{0}> \xi_{1} >  \dots > \xi_{2k} <\xi_{2k+1}  \right)\\
=  \int_{0}^t(1-F(r)) dF(r) + \sum_{k\geq 1} \int_{0}^t \mathbb{P}\left(t\geq \xi_{0}> \xi_{1} >  \dots > \xi_{2k-1}> r\ ;\ r <\xi_{2k+1}\right) dF(r),\\
=  \int_{0}^t(1-F(r)) dF(r) + \sum_{k\geq 1} \int_{0}^t \mathbb{P}\left(t\geq \xi_{0}> \xi_{1} >  \dots > \xi_{2k-1}> r\right)\mathbb{P}\left(r <\xi_{2k+1}\right) dF(r),
\end{multline*}
at second line we have conditioned respectively on $\set{\xi_0=r}$ and on $\set{\xi_{2k}=r}$. 

Set $A=\set{\xi_{0}, \xi_{1},\dots, \xi_{2k-1}\in (r,t)}$, then $\mathbb{P}(A)=(F(t)-F(r))^{2k}$ and conditional on $A$ these random variables are ordered as a uniform permutation:
\begin{align*}
f(t) &=  \int_{r=0}^t (1-F(r)) dF(r) + \sum_{k\geq 1} \int_{0}^t \frac{1}{(2k)!} (F(t)-F(r))^{2k} (1-F(r)) dF(r)\\
&=\int_{0}^t  \sum_{k\geq 0} \frac{1}{(2k)!}   (F(t)-F(r))^{2k} (1-F(r)) dF(r)\\
&=\int_{0}^{F(t)} \sum_{k\geq 0} \frac{1}{(2k)!}  (F(t)-s)^{2k} (1-s) ds\\
&=\int_{0}^{F(t)}\cosh\left(F(t)-s\right) (1-s)ds\\
&=1-\exp\big(-F(t)\big).
\end{align*}
\end{proof}

\section{Parking on an interval}
\subsection{Coupling with the infinite parking}
Recall that $M_n$ is the number of places occupied in the finite parking of size $n$.
\begin{theorem}
For every $n\geq 2$,
\begin{equation}\label{Eq:LimiteVoituresRaffinee}
\left|\mathbb{E}[M_{n}]-n(1-e^{-2})\right|\leq 14.
\end{equation}
\end{theorem}
This estimate is not as sharp as $\mathbb{E}[M_n]=n(1-e^{-2})+ (1-3e^{-2})+\mathrm{o}(1)$ which has been proved by Friedman \cite{Friedman} (see also \cite{Flajolet}). Yet the proof we provide here provides a simple estimate on $\mathbb{P}(X_n(i)=0)$ which will be useful later.
We use $\xi_i$'s to define a \emph{coupling} between finite and infinite parkings. 

\begin{proof} 
Let $(\xi_i)_{i\in\mathbb{Z}}$ be as before a sequence of i.i.d. continuous random variables of common distribution function $F$, we use the same $\xi_i$'s to define $(X_n(i))_{1\leq i\leq n}$ and  $(X_{\infty}(i))_{i\in\mathbb{Z}}$. Let $m$ (resp. $m'$) be the leftmost (resp. rightmost) local minimum of $(\xi_i)_{i\in\mathbb{Z}}$ in $\set{1,\dots,n-1}$. We set $m=n$, $m'=0$ if there is no local minimum.

By construction, if there is a local minimum, $X_n(i)$ and $X_{\infty}(i)$ coincide for every $m\leq i\leq m'+1$. Therefore
\begin{align*}
\left|\mathbb{E}[M_{n}] -n(1-e^{-2})\right|&= \left|\sum_{i=1}^n \mathbb{E}[X_n(i)]-\mathbb{E}[X_\infty(i)]\right|\\
&\leq \mathbb{E}\left[\mathrm{card}\set{i;\ X_n(i)\neq X_{\infty}(i)}\right]\\
&\leq \mathbb{E}[m-1]+\mathbb{E}[n-m'-1]=2\mathbb{E}[m-1].
\end{align*}
There are different ways to bound $\mathbb{E}[m]$, an easy one is to observe that local minima appear independently at $2,5,8,\dots,$ with probability $\mathbb{P}(\xi_1>\xi_2<\xi_3)=1/3$.
We get $\mathbb{P}\left(m\geq i\right)\leq (2/3)^{\lf i/3\rf}$ which yields $\mathbb{E}[m-1]\leq 7$ for every $n$. 
\end{proof}
Note that as a by-product of the proof we get that
\begin{equation}
\left|\mathbb{P}(X_{n}(i)=0)-e^{-2}\right|\leq \mathbb{P}\left(m>i\ \cup\ m'<i-1 \right)
\leq 2 \eps_i,\label{Eq:MajoExplicite2}
\end{equation}
where $\eps_i=\max\set{ (2/3)^{i/3-1},(2/3)^{(n-i)/3-1} }$.

We also can deduce the convergence in probability of $M_n/n$ by using
$$
\sum_{i=1}^n |X_n(i)-X_\infty(i)| \leq m-1+(n-m'-1)
$$
and the fact that $m/n,(n-m')/n$ both converges to zero in probability.

\subsection{Number of trials: Poissonization}\label{Sec:Coupons}
Let $T_n$ be the number of cars that have tried to park before the parking process is over. It is clear that $T_n$ is stochastically smaller than the number of trials needed to pick each number in $\set{1,\dots, n-1}$ at least once, \emph{i.e.} stochastically smaller than a coupon collector with $n-1$ coupons. Thus, the $\limsup$ (in probability) of $T_n/(n\log n)$ is less than one. 

In order to estimate $T_n$ we use another construction of the arrival process, in order to take into account the arrivals of cars that tried but did not succeed in park. We now are given for each $i\in \set{1,\dots,n}$ a sequence of random variables $(\xi^j_i)_{j\geq 1}$, the family $\big\{\xi^j_i\big\}_{i,j}$ being i.i.d. exponentially distributed with mean one. At $(i,i+1)$, cars try to park at times
$$
\xi^1_i,\xi^1_i+\xi^2_i, \xi^1_i+\xi^2_i+\xi^3_i,\dots
$$
For simplicity, we write as before $\xi_i=\xi_i^1$ for the first arrival of a car at $i$.
Let $\tau_\star$ be the arrival of the last car that succeeds in parking:
$$
\tau_\star=\max \set{\tau_i, \tau_i<+\infty}.
$$
Here is a picture that sums up notations (here the last car parks at $(i,i+1)$, we have $T_n=12$):

\medskip

\begin{center}
\includegraphics[width=12cm]{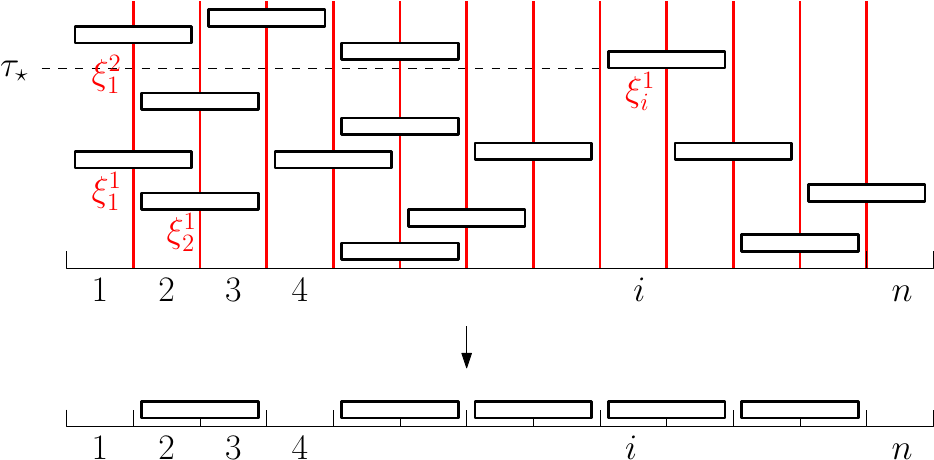}\\
\end{center}

\medskip

By construction and by Markovianity  we have
$$
T_n \stackrel{(d)}{=} \sum_{i=1}^n \mathrm{max}\set{j;\ \xi^1_i+\dots +\xi^j_i \leq \tau_\star},
$$
since $\mathrm{max}\set{j;\ \xi^1_i+\dots +\xi^j_i \leq t}$ is the number of cars that tried to park at $(i,i+1)$ before time $t$.

\begin{theorem}[Number of trials]
$$
\frac{T_n}{n\log(n)} \stackrel{\text{prob.}}{\to} 1.
$$
\end{theorem}
\begin{proof}

\noindent{\bf Upper bound.} 
As noted above, $T_n$ is stochastically smaller than a coupon collector with $n-1$ coupons. It is classical (see for instance \cite{Durrett} Example 2.2.3) that for each $\eps>0$ we have
$$
\mathbb{P}\left(\frac{T_n}{n\log(n)} \geq 1+\eps \right) \to 0
$$
\noindent{\bf Lower bound.} 
The strategy is the following: Theorem \ref{Th:EvolutionDensite} suggests that, as long as $i$ is bounded away from $0$ and $n$, we have $\mathbb{P}(u\leq \tau_i<+\infty)\approx e^{-2(1-e^{-u})}-e^{-2} \sim_{u\to +\infty} 2e^{-2}e^{-u}$ . The $\tau_i$'s being weakly dependent, we expect $\tau_\star =\max\set{\tau_i;\ \tau_i <+\infty}$ to be of order $\log(n)$.
To conclude, we will use the fact that $T_n\approx n\times \tau_\star$. 

\begin{lemma}\label{Lem:EstimeeTauStar}
For every $\delta >0$,
$$
\mathbb{P}(\tau_\star \geq (1-\delta)\log(n)) \stackrel{n\to +\infty}{\to} 1.
$$
\end{lemma}
\begin{proof}[Proof of Lemma \ref{Lem:EstimeeTauStar}]
For two integers $i,\ell$ such that $[i-\ell,i+\ell]\subset [1,n]$, let $A_{i,\ell}(u)$ be the event $E_{i,\ell}(u)\cup F_{i,\ell}(u)\cup F_{i,\ell}(u) $ where
\begin{align*}
E_{i,\ell}(u)= \big\{&\text{ odd rise of length $\leq \ell-2$ at $i$},\ \xi_{i-1} \geq u,\\ &\text{ even descent of length $\leq \ell-2$ at $i$}\ \big\}\\
F_{i,\ell}(u)= \big\{&\text{ even rise of length $\leq \ell-2$ at $i$},\ \xi_{i} \geq u,\\ &\text{ odd descent of length $\leq \ell-2$ at $i$}\ \big\}\\
G_{i,\ell}(u)= \big\{&\text{ odd rise of length $\leq \ell-2$ at $i$},\ \xi_{i-1} \geq u,\\ &\xi_i\geq u,\text{ odd descent of length $\leq \ell-2$ at $i$}\ \big\}.
\end{align*}
The intuition for $A_{i,\ell}(u)$ is that, at least for large $\ell$, this event says that a car will arrive at $i$, but not before time $u$.

Event $A_{i,\ell}$ only depends on $\set{\xi_{i'}, i-\ell+1\leq i'\leq i+\ell-1}$.
Again by Lemma \ref{Lem:RiseDescent} we have
$$
\set{u\leq \tau_i<+\infty}
\supset
A_{i,\ell}(u)
\supset
\begin{cases}
u\leq \tau_i<+\infty,\\
\text{{\bf and} there is a local minimum among $\xi_{i-\ell+2},\dots,\xi_{i-1}$}\\
\text{{\bf and} there is a local minimum among $\xi_i,\dots,\xi_{i+\ell-2}$}.
\end{cases}
$$
Then
\begin{align*}
0\leq  \mathbb{P}(u\leq \tau_i<+\infty )-\mathbb{P}(A_{i,\ell}(u))&\leq \mathbb{P}\bigg(\text{ no local min. between $i-\ell+2$ and $i-1$ }\\
&\phantom{222222}\text{{\bf or} no local min. between $i$ and $i+\ell-2$}\bigg)\\
&\leq 2(2/3)^{\lf (\ell-2)/3\rf} \leq 2(2/3)^{ \ell/4}
\end{align*}
for large $\ell$. (Here we have re-used the fact that local minima appear independently at $i+1,i+4,i+7,\dots,i+\ell$ with probability $1/3$.)

Besides,
\begin{align}
\mathbb{P}(u \leq \tau_i <+\infty)=\ &1-\mathbb{P}(\tau_i=+\infty)-\mathbb{P}(\tau_i< u)\notag\\
=\ &1-\mathbb{P}(\tau_i=+\infty)\\
&-\mathbb{P}\bigg( \big\{\xi_{i-1}< u; \text{ odd rise at $i$ }\big\}\cup \big\{\xi_{i}< u; \text{ odd descent at $i$ }\big\}\bigg)\label{Eq:GrosseEquation}
\end{align}
Aside from the boundary effects, the arguments of the proof of Theorem \ref{Th:EvolutionDensite} are still valid and we get for $i\geq 2$
\begin{align*}
\mathbb{P}\big(  \xi_{i-1} < u; \text{ odd rise at $i$ }\big)
&=\sum_{k=0}^{\lf i/2\rf -1} \mathbb{P}\left(\xi_{i-(2k+2)}>\xi_{i-(2k+1)}< \xi_{i-2k} <  \dots < \xi_{i-1} \leq u \right)\\
&=\int_{0}^{F(u)} \sum_{k= 0}^{\lf i/2\rf -1} \frac{1}{(2k)!}  (s-F(t))^{2k} (1-s) ds\\
&=\delta_i+\int_{0}^{F(u)}\cosh\left(s-F(t)\right) (1-s)ds \\
&=1-\exp\big(-F(u)\big)+\delta_i,
\end{align*}
where $|\delta_i|\leq 2/i!$ ($\delta_i$ is obtained by bounding the remainder of the Taylor series of $\cosh$). By symmetry $i\leftrightarrow n-i$ we have the symmetric estimate on $\mathbb{P}\big(  \xi_{i}\leq u; \text{ odd descent at $i$ }\big)$.

Plugging this into \eqref{Eq:GrosseEquation} and combining with our estimate \eqref{Eq:MajoExplicite2} on $\mathbb{P}(\tau_i=+\infty)$  we obtain for $i\geq 2$, since $\eps_i \geq \delta_i$
$$
\mathbb{P}(u\leq \tau_i <+\infty )\geq e^{-2F(u)}-e^{-2} -10\eta_i\geq 2e^{-2}(1-F(u))-10\eta_i= 2e^{-2}e^{-u}-10\eta_i,
$$
where $|\eta_i| \leq \max\set{ (2/3)^{i/3-1},(2/3)^{(n-i)/3-1} }$.

Now, events
$$
A_{\ell,\ell}(u),A_{3\ell,\ell}(u),A_{5\ell,\ell}(u),\dots A_{\lf n/\ell-1\rf \ell,\ell}(u)
$$
are independent and (we skip integer parts in order to lighten notations):
\begin{align*}
\mathbb{P}(\tau_\star \leq u)&\leq \mathbb{P}\bigg(\text{not }A_{\ell,\ell}(u),\text{not }A_{3\ell,\ell}(u),\text{not }A_{5\ell,\ell}(u),\dots , \text{not }A_{n/\ell\times \ell,\ell}(u)\bigg)\\
&\leq \prod_{j=1}^{ n/\ell} \left(1-2e^{-2}e^{-u}+10\eta_j +2(2/3)^{\ell/4}\right) \\
&\leq \prod_{j=\log(n)}^{ n/\ell-\log(n)} \left(1-2e^{-2}e^{-u}+10\eta_j +2(2/3)^{\ell/4}\right).
\end{align*}
Choose now $\ell=50\log(n)$ and take $u=(1-\delta)\log(n)$, so that for large $n$ every term of the product is less than $1-e^{-2}e^{-u}$,
$$
\mathbb{P}(\tau_\star \leq (1-\delta)\log(n))\leq
\left(1-\tfrac{e^{-2}}{n^{1-\delta}}\right)^{\left(n/50\log(n)-2\log(n)\right)}\leq \exp(-n^{\delta/2}),
$$
for large $n$.
\end{proof}
We now conclude the lower bound:
\begin{align*}
\mathbb{P}(T_n\leq (1-\eps)n&\log(n)) \\
\leq\ &\mathbb{P}(\tau_\star \leq (1-\delta)\log(n))\\
&+ \mathbb{P}(T_n\leq (1-\eps)n\log(n) ; \tau_\star > (1-\delta)\log(n))\\
\leq\ &\mathbb{P}(\tau_\star \leq (1-\delta)\log(n)) \\
&+\mathbb{P}\left( \sum_{i=1}^n \max\set{j;\ \xi^1_i+\dots +\xi^j_i \leq (1-\delta)\log(n)} \leq (1-\eps)n\log(n)\right)\\
\leq\ &\mathbb{P}(\tau_\star \leq (1-\delta)\log(n)) \\
&+\mathbb{P}\left( \sum_{i=1}^n \mathrm{Poiss}_i\left((1-\delta)\log(n)\right) \leq (1-\eps)n\log(n)\right),
\end{align*}
where $\mathrm{Poiss}_i(\lambda)$ are i.i.d. Poisson with mean $\lambda$. The first term in the right-hand side goes to zero thanks to Lemma \ref{Lem:EstimeeTauStar}, so does the second one by taking $\delta=\eps/2$ and using Chebyshev's inequality.
\end{proof}

\subsection*{Acknowledgements}
I learned about the Page model in the very nice book \cite{Pinsky} by Ross Pinsky. I also would like to thank Marie Albenque and Jean-Fran\c{c}ois Marckert for stimulating discussions, and an anonymous referee for her/his careful reading.

\end{document}